\documentclass[12pt]{amsart}
\usepackage{a4wide}
\usepackage[colorlinks]{hyperref}
\AtBeginDocument{
  \hypersetup{
    linkcolor=blue,
    citecolor=green,
  }
}
\usepackage{mathscinet}
\usepackage{mathtools}
\usepackage{bm}
\usepackage[T1]{fontenc}

\makeatletter
\let\c@author\relax
\makeatother
\usepackage[backend=bibtex,style=numeric,sorting=nyt,
giveninits=true,isbn=false,url=false,maxbibnames=99]{biblatex}
\renewbibmacro{in:}{\ifentrytype{article}{}{\printtext{\bibstring{in}\intitlepunct}}}
\DeclareFieldFormat*{title}{#1}
\usepackage[mathlines,pagewise]{lineno}
\newcommand*\patchAmsMathEnvironmentForLineno[1]{%
  \expandafter\let\csname old#1\expandafter\endcsname\csname #1\endcsname
  \expandafter\let\csname oldend#1\expandafter\endcsname\csname end#1\endcsname
  \renewenvironment{#1}%
     {\linenomath\csname old#1\endcsname}%
     {\csname oldend#1\endcsname\endlinenomath}}%
\newcommand*\patchBothAmsMathEnvironmentsForLineno[1]{%
  \patchAmsMathEnvironmentForLineno{#1}%
  \patchAmsMathEnvironmentForLineno{#1*}}%
\AtBeginDocument{%
\patchBothAmsMathEnvironmentsForLineno{equation}%
\patchBothAmsMathEnvironmentsForLineno{align}%
\patchBothAmsMathEnvironmentsForLineno{flalign}%
\patchBothAmsMathEnvironmentsForLineno{alignat}%
\patchBothAmsMathEnvironmentsForLineno{gather}%
\patchBothAmsMathEnvironmentsForLineno{multline}%
}

\usepackage{cleveref}
\crefname{equation}{}{}

\newtheorem{theorem}{Theorem}[section]
\newtheorem{lemma}[theorem]{Lemma}

\newtheorem{corollary}[theorem]{Corollary}

\theoremstyle{definition}

\theoremstyle{remark}
\newtheorem{remark}[theorem]{Remark}

\numberwithin{equation}{section}

\begin{document}

\title[Interior $C^{2,\alpha}$ regularity]{Interior $C^{2,\alpha}$ regularity for fully nonlinear uniformly elliptic equations in dimension two}

\author{Kai Zhang}
\address{Departamento de Geometr\'{i}a y Topolog\'{i}a,
Instituto de Matem\'{a}ticas IMAG, Universidad de Granada}
\email{zhangkaizfz@gmail.com; zhangkai@ugr.es;\newline
 MR ID:\href{https://mathscinet.ams.org/mathscinet/2006/mathscinet/search/author.html?mrauthid=1098004}
{1098004};
ORCID:\href{https://orcid.org/0000-0002-1896-3206}
{0000-0002-1896-3206}
}

\thanks{This research has been financially supported by the Project PID2020-118137GB-I00 funded by MCIN/AEI /10.13039/501100011033.}


\date{}



\begin{abstract}
In this note, we present the interior $C^{2,\alpha}$ regularity for viscosity solutions of fully nonlinear uniformly elliptic equations in dimension two.
\end{abstract}

\maketitle

\section{Introduction}\label{S1}
In this note, we give the interior $C^{2,\alpha}$ regularity for viscosity solutions of the following fully nonlinear uniformly elliptic equation:
\begin{equation}\label{FNE}
F(D^2u,Du,u,x)=f\quad\mbox{in}~~B_1,
\end{equation}
where $B_1\subset \mathbb{R}^2$ is the unit disk. The main assumption on $F$ is the uniform ellipticity, i.e., $F$ is continuous in all its arguments and for any $M,N\in \mathcal{S}^2$ (the set of $2\times 2$ symmetric matrices), $p,q\in \mathbb{R}^2$, $s,t\in \mathbb{R}$ and $x\in B_1$,
\begin{equation}\label{SC1}
  \begin{aligned}
F(M,p,s,x)-F(N,q,t,x)\leq \Lambda |(M-N)^+|-\lambda |(M-N)^-|+b|p-q|+c|s-t|,
  \end{aligned}
\end{equation}
where $0<\lambda\leq \Lambda$ and $b,c\geq 0$ are constants; the $(M-N)^+$ and $(M-N)^-$ denote the positive part and negative part of $M-N$ respectively.

According to Evans-Krylov theory (see \cite{MR649348}, \cite{MR661144}, or \cite[Theorem 6.6]{MR1351007}), interior $C^{2,\alpha}$ regularity holds provided $F$ is convex or concave in $M$. However, it is also well known that this assumption is unnecessary in dimension two, a result established by Nirenberg \cite[Theorem I]{MR0064986}. In \cite[Theorem I]{MR0064986}, it requires that $F\in C^1$ and $u\in C^2$. The primary objective of this note is to remove these smoothness assumptions. This is possible because the $C^{2,\alpha}$ estimate is independent of these smoothness conditions. Although this result may be well-known to experts in the field, we have been unable to locate an explicit reference in the existing literature. Therefore, we present the proof here for the sake of completeness.

The main result in this note is:
\begin{theorem}\label{th1.1}
Suppose that $u\in C(\bar{B}_1)$ is a viscosity solution of
\begin{equation*}
F(D^2u)=0\quad\mbox{in}~~B_{1}\subset \mathbb{R}^2,
\end{equation*}
where $F$ is uniformly elliptic and $F(0)=0$. Then there exists a universal constant $0<\bar \alpha<1$ such that $u\in C^{2,\bar \alpha}(\bar{B}_{1/2})$ and
\begin{equation}\label{e1.1}
\|u\|_{C^{2,\bar \alpha}(\bar{B}_{1/2})} \leq C\|u\|_{L^{\infty }(B_1)},\\
\end{equation}
where $C$ is universal. A constant is called universal if it depends only on $\lambda$ and $\Lambda$.
\end{theorem}

Once we have \Cref{th1.1}, we can derive the following $C^{2,\alpha}$ regularity for general equations by the perturbation argument (see \cite[Theorem 5.1, Remark 5.2 and Corollary 8.4]{lian2020pointwise}).
\begin{corollary}\label{co1.1}
Let $u\in C(\bar B_1)$ be a viscosity solution of \cref{FNE} where $F$ is uniformly elliptic. Suppose that for some $0<\alpha<\bar \alpha$, $f\in C^{\alpha}(\bar{B}_1)$ and
\begin{equation*}
\begin{aligned}
|F(M,p,s,x)-F(M,p,s,y)|\leq &K_1|x-y|^{\alpha} \left(|M|+|p|+|s|\right)+K_2|x-y|^{\alpha} ,\\
&~\forall ~(M,p,s)\in \mathcal{S}^2\times \mathbb{R}^2 \times \mathbb{R},~~ x,y\in B_1.
\end{aligned}
\end{equation*}

Then $u\in C^{2,\alpha}(\bar{B}_{1/2})$ and
\begin{equation*}
\|u\|_{C^{2,\alpha}(\bar{B}_{1/2})}\leq C\left(\|u\|_{L^{\infty}(B_1)}+\|f\|_{C^{\alpha}(\bar{B}_{1})}+|F(0,0,0,0)|+K_2\right),
\end{equation*}
where $C$ depends only on $\lambda,\Lambda,b,c,K_1$ and $\alpha$.
\end{corollary}

\section{A priori \texorpdfstring{$C^{2,\alpha}$}{C2,a} estimate}
In this section, we obtain the global a priori $C^{2,\alpha}$ estimate under smooth assumptions. The following is the main result in this section:
\begin{theorem}\label{th3.3}
Let $u\in C^2(\bar{B}_1)$ be a classical solution of
\begin{equation*}
\left\{\begin{aligned}
&F(D^2u)=0&& ~~\mbox{in}~~B_1;\\
&u=g&& ~~\mbox{on}~~\partial B_1.
\end{aligned}\right.
\end{equation*}
Suppose that $F\in C^1$ is uniformly elliptic, $F(0)=0$ and $g\in C^{2,\bar \alpha}(\partial B_1)$. Then $u\in C^{2,\bar \alpha}(\bar{B}_1)$ and
\begin{equation}\label{e3.2}
\|u\|_{C^{2,\bar \alpha}(\bar{B}_1)}\leq C\|g\|_{C^{2,\bar \alpha}(\partial B_1)},
\end{equation}
where $0<\bar \alpha<1$ and $C$ are universal.
\end{theorem}

We start from the cornerstone of this note: interior $C^{2,\alpha}$ estimate for smooth solutions, which is the due to Nirenberg \cite[Theorem I]{MR0064986} (see also \cite[Theorem 4.9]{MR4560756}).
\begin{theorem}\label{th3.1}
There exists a universal constant $0<\bar \alpha<1$ such that the following holds.  Let $u\in C^2(\bar{B}_1)$ be a classical solution of
\begin{equation}\label{e3.1}
F(D^2u)=0\quad\mbox{in}~~B_{1},
\end{equation}
where $F\in C^1$ is uniformly elliptic and $F(0)=0$. Then $u\in C^{2,\bar \alpha}(\bar{B}_{1/2})$ and
\begin{equation}\label{e3.1-2}
\|u\|_{C^{2,\bar \alpha}(\bar{B}_{1/2})}
\leq C\|u\|_{L^{\infty }(B_1)},
\end{equation}
where  $C$ is universal.
\end{theorem}

The next is the boundary $C^{2,\alpha}$ regularity for fully nonlinear uniformly elliptic equations (see \cite[Theorem 1.2]{MR3246039} or \cite[Theorem 1.8]{MR4088470})). We remark here that this result holds for any dimension without the convexity/concavity assumption on $F$.
\begin{theorem}\label{th3.2}
There exists a universal constant $0<\bar \alpha<1$ such that the following holds. Let $u$ be a viscosity solution of
\begin{equation*}
\left\{\begin{aligned}
&F(D^2u)=0&& ~~\mbox{in}~~B_1\subset \mathbb{R}^2;\\
&u=g&& ~~\mbox{on}~~\partial B_1
\end{aligned}\right.
\end{equation*}
and $x_0\in \partial B_1$. Suppose that $F$ is uniformly elliptic, $F(0)=0$ and $g\in C^{2,\bar \alpha}(x_0)$.

Then $u\in C^{2,\bar \alpha}(x_0)$, i.e., there exists a quadratic polynomial $P_{x_0}$ such that
\begin{equation}\label{en2}
  |u(x)-P_{x_0}(x)|\leq C |x-x_0|^{2+\bar \alpha}\left(\|u\|_{L^{\infty }(B_1)}+\|g\|_{C^{2,\bar \alpha}(x_0)}\right), ~~\forall ~x\in B_1,
\end{equation}
\begin{equation}\label{e2.0}
F(D^2P_{x_0})=0
\end{equation}
and
\begin{equation}\label{e.t2-2}
|DP_{x_0}(x_0)|+|D^2P_{x_0}|\leq C\left(\|u\|_{L^{\infty }(B_1)}+\|g\|_{C^{2,\bar \alpha}(x_0)}\right),
\end{equation}
where $C$ is universal.
\end{theorem}

\begin{remark}\label{re2.2}
In the theorem above, we do not assume $F \in C^1$, unlike in \Cref{th3.1}. Note that the H\"{o}lder exponent $\bar \alpha$ in \Cref{th3.2} may differ from that in \Cref{th3.1}. However, by choosing the smaller of the two, we ensure that both theorems hold for a common $\bar \alpha$. From now on, we fix such an $\bar \alpha$ so that both \Cref{th3.1} and \Cref{th3.2} remain valid.
\end{remark}

It is standard that by combining \Cref{th3.1} and \Cref{th3.2}, we have the global $C^{2,\alpha}$ estimates:~\\
\noindent\textbf{Proof of \Cref{th3.3}.} Without loss of generality, we assume that
\begin{equation*}
\|g\|_{C^{2,\bar \alpha}(\partial B_1)}\leq 1.
\end{equation*}
By the maximum principle,
\begin{equation*}
\|u\|_{L^{\infty}(B_1)}\leq \|g\|_{L^{\infty}(\partial B_1)}\leq 1.
\end{equation*}

To prove \cref{e3.2}, we only need to prove that given any $x_0\in \bar{B}_1$, there exists a quadratic polynomial $P_{x_0}$ such that
\begin{equation}\label{e.l61-3}
  |u(x)-P_{x_0}(x)|\leq C|x-x_0|^{2+\bar \alpha},~\forall ~x\in \bar{B}_1,
\end{equation}
and
\begin{equation}\label{e2.1}
|DP_{x_0}(x_0)|+|D^2P_{x_0}|\leq C
\end{equation}
where $C$ is universal. Throughout this proof, $C$ always denotes a universal constant.

If $x_0\in \partial B_1$, \cref{e.l61-3} follows from \Cref{th3.2} directly. In the following, we assume $x_0\in B_1$.
Let $\tilde{x}_0\in \partial B_1$ with
\begin{equation*}
r_0:=|x_0-\tilde{x}_0|=d(x_0,\partial B_1),
\end{equation*}
where $d(x_0,\partial B_1)$ denotes the distance from $x_0$ to $\partial B_1$.

By the boundary $C^{2,\alpha}$ regularity \Cref{th3.2}, there exists $P_{\tilde{x}_0}$ such that
\begin{equation}\label{e9.4}
  |u(x)-P_{\tilde{x}_0}(x)|\leq C|x-\tilde{x}_0|^{2+\bar \alpha}, ~\forall ~x\in \bar{B}_{1}.
\end{equation}
Set
\begin{equation*}
v=u-P_{\tilde{x}_0}
\end{equation*}
and then $v$ satisfies
\begin{equation}\label{e.Cka-4}
  \tilde{F}(D^2v)=0\quad\mathrm{in}~B_{r_0}(x_0),
\end{equation}
where
\begin{equation*}
\tilde{F}(M)\coloneqq F(M+D^2P_{\tilde{x}_0}).
\end{equation*}
Moreover, $\tilde{F}(0)=F(D^2P_{\tilde{x}_0})=0$ (see \Cref{e2.0})

Then \cref{e.Cka-4} satisfies the assumptions of the interior $C^{2,\alpha}$ regularity (see \Cref{th3.1}). Hence, there exists a quadratic polynomial $P$ such that (with the aid of \cref{e9.4} and the definition of $v$)
\begin{equation}\label{e9.5}
  \begin{aligned}
|v(x)-P(x)|\leq& C\|v\|_{L^\infty(B_{r_0}(x_0))}r_0^{-(2+\bar \alpha)}|x-x_0|^{2+\bar \alpha}
  \leq C |x-x_0|^{2+\bar \alpha}\quad\mathrm{in}~B_{r_0/2}(x_0),\\
|P(x_0)|=&|v(x_0)|=|u(x_0)-P_{\tilde{x}_0}(x_0)|\leq C|x_0-\tilde{x}_0|=Cr_0^{2+\bar \alpha},\\
|D P(x_0)|\leq& Cr_0^{-1}\|v\|_{L^\infty(B_{r_0}(x_0))}\leq Cr_0^{1+\bar \alpha},\\
|D^2P(x_0)|\leq& Cr_0^{-2}\|v\|_{L^\infty(B_{r_0}(x))}\leq Cr_0^{\bar \alpha}.
  \end{aligned}
\end{equation}

Let
\begin{equation*}
P_{x_0}=P_{\tilde{x}_0}+P.
\end{equation*}
For any $x\in \bar{B}_1$, if $|x-x_0|<r_0/2$, by the first inequality in \cref{e9.5},
\begin{equation*}
  |u(x)-P_{x_0}(x)|=|v(x)-P(x)|\leq C|x-x_0|^{2+\bar \alpha}.
\end{equation*}
If $|x-x_0|\geq r_0/2$, we have
\begin{equation*}
  \begin{aligned}
|u(x)-P_{x_0}(x)|&\leq |u(x)-P_{\tilde{x}_0}(x)|+|P(x)|\\
&\leq C|x-\tilde{x}_0|^{2+\bar \alpha}
+\left(|P(x_0)|+|DP(x_0)|\cdot|x-x_0|+|D^2P(x_0)|\cdot|x-x_0|^2\right)\\
&\leq C|x-\tilde{x}_0|^{2+\bar \alpha}+C(r_0^{2+\bar \alpha}+|x-x_0|r_0^{1+\bar \alpha}+|x-x_0|^2r_0^{\bar \alpha})\\
&\leq C|x-x_0|^{2+\bar \alpha}.
  \end{aligned}
\end{equation*}
Therefore, \cref{e.l61-3} holds for any $x\in \bar{B}_1$.~\qed

\section{Existence of classical solutions}\label{S3}
In this section, we prove the existence of classical solutions by the method of continuity with the aid of a priori $C^{2,\alpha}$ estimate obtained in last section. The main result is the following.
\begin{theorem}\label{le4.1}
Let $F\in C^2$ be uniformly elliptic with $F(0)=0$. Suppose that $g\in C^{2,\bar \alpha}(\partial B_1)$. Then there exists a unique classical solution $u\in C^{2,\bar \alpha}(\bar B_1)$ of
\begin{equation}\label{e4.0}
\left\{\begin{aligned}
&F(D^2u)=0&& ~~\mbox{in}~~B_1\subset \mathbb{R}^2;\\
&u=g&& ~~\mbox{on}~~\partial B_1.
\end{aligned}\right.
\end{equation}
Moreover, we have
\begin{equation}\label{e3.3}
\|u\|_{C^{2,\bar \alpha}(\bar{B}_{1})}\leq C\|g\|_{C^{2,\bar \alpha}(\partial B_{1})},
\end{equation}
where $C$ is universal.
\end{theorem}
\proof The uniqueness follows from the maximum principle (see \cite[Corollary 5.4]{MR1351007}) and the estimate \cref{e3.3} follows from \Cref{th3.3}. Hence, we only need to prove the existence. We use the method of continuity (see \cite[Proof of Theorem 9.7]{MR1351007}). We consider a family of Dirichlet problems for $t\in [0,1]$:
\begin{equation}\label{e4.1}
\left\{\begin{aligned}
&tF(D^2u)+(1-t)\Delta u=0&& ~~\mbox{in}~~B_1;\\
&u=g&& ~~\mbox{on}~~\partial B_1.
\end{aligned}\right.
\end{equation}
Let
\begin{equation*}
A:=\left\{t\in [0,1]: ~~\mbox{there exists a solution}~~ u\in C^{2,\bar \alpha}(\bar{B}_1)~~\mbox{of}~~\cref{e4.1}\right\}.
\end{equation*}
Clearly, $0\in A$. If we can show that $A$ is both open and closed, we have $1\in A$, which means that there exists a solution of \cref{e4.0}.

We first show that $A$ is open. Suppose that $t_0\in A$. Then there exists a solution $u_0$ of \cref{e4.1}. Let
\begin{equation*}
C^{2,\bar \alpha}_0(\bar{B}_1):=\left\{u\in C^{2,\bar \alpha}(\bar{B}_1):~~u=0\quad\mbox{on}~~\partial B_1\right\}.
\end{equation*}
Define an operator $\phi: C^{2,\bar \alpha}_0(\bar{B}_1)\times \mathbb{R}\to C^{\bar \alpha}(\bar{B}_1)$ as follows:
\begin{equation*}
\phi(u,t):=tF(D^2u_0+D^2u)+(1-t)\Delta (u_0+u).
\end{equation*}
Clearly, $\phi(0,t_0)=0$. Now, we prove that $D_u\phi (0,t_0):C^{2,\bar \alpha}_0(\bar{B}_1)\to C^{\bar \alpha}(\bar{B}_1)$ is an invertible bounded linear operator. Indeed, by a direct calculation,
\begin{equation}\label{e4.2}
\langle D_u\phi (0,t_0), w\rangle=tF_{M_{ij}}(D^2u_0)w_{ij}+(1-t)\Delta w,~\forall ~w\in C^{2,\alpha}_0(\bar{B}_1).
\end{equation}
Let
\begin{equation*}
a^{ij}(x):= tF_{M_{ij}}(D^2u_0(x))+(1-t)\delta^{ij},~\forall ~x\in B_1,
\end{equation*}
where $\delta^{ij}$ is the Kronecker symbol and denotes the unit matrix. Then \cref{e4.2} can be rewritten as
\begin{equation*}
\langle D_u\phi (0,t_0), w\rangle=a^{ij}w_{ij},
\end{equation*}
where the Einstein summation convention is used.

Since $F\in C^2$ and $u_0\in C^{2,\bar \alpha}(\bar{B}_1)$, we have $a^{ij}\in C^{\bar \alpha}(\bar{B}_1)$. By the classical Schauder theory for linear uniformly elliptic equation (see \cite[Theorem 6.14]{MR1814364}), for any $f\in C^{\bar \alpha}(\bar{B}_1), $ there exists a unique solution $u\in C^{2,\bar \alpha}(\bar{B}_1)$ of
\begin{equation*}
\left\{\begin{aligned}
&a^{ij}w_{ij}=f&& ~~\mbox{in}~~B_1;\\
&u=0&& ~~\mbox{on}~~\partial B_1
\end{aligned}\right.
\end{equation*}
and
\begin{equation*}
\|u\|_{C^{2,\bar \alpha}(\bar{B}_{1})}\leq C\|f\|_{C^{\bar \alpha}(\bar B_{1})},
\end{equation*}
where $C$ depends on $\lambda,\Lambda, u_0$ and the smoothness of $F$. Thus, $D_u\phi (0,t_0)$ is an invertible linear operator from $C^{2,\bar \alpha}_0(\bar{B}_1)$ to $C^{\bar \alpha}(\bar{B}_1)$. By the implicit function theorem, there exists a neighborhood $(t_0-\varepsilon,t_0+\varepsilon)$ of $t_0$ such that for any $t\in (t_0-\varepsilon,t_0+\varepsilon)$, there exists a unique $u\in C^{2,\bar \alpha}_0(\bar{B}_1)$ such that $\phi(u,t)=0$. In another word, $u_0+u$ is a solution of \cref{e4.1}. Therefore, we have proved that $A$ is open, i.e., if $t_0\in A$, there exists a neighborhood $(t_0-\varepsilon,t_0+\varepsilon)\cap [0,1]\subset A$.

Next, we prove that $A$ is closed. Given a sequence of $t_m\in A$ with $t_m\to \bar t$, we need to prove that $\bar t\in A$. By the definition of $A$, there exists a sequence of solutions $u_m$ of \cref{e4.1} with $t=t_m$ there. By the $C^{2,\alpha}$ estimate (see \Cref{th3.3}),
\begin{equation*}
\|u_m\|_{C^{2,\bar \alpha}(\bar{B}_1)}\leq C\|g\|_{C^{2,\bar \alpha}(\partial B_1)},
\end{equation*}
where $C$ is independent of $m$. Since $C^{2,\bar \alpha}(\bar{B}_1)$ is compactly embedded into $C^{2}(\bar{B}_1)$, there exists a subsequence (denoted by $u_m$ again) and $\bar{u}\in C^{2,\bar \alpha}(\bar{B}_1)$ such that
\begin{equation*}
u_m\to \bar{u}\quad\mbox{in}~~C^{2}(\bar{B}_1)
\end{equation*}
and
\begin{equation*}
\|\bar{u}\|_{C^{2,\bar \alpha}(\bar{B}_1)}\leq C\|g\|_{C^{2,\bar \alpha}(\partial B_1)}.
\end{equation*}

Take the limit in
\begin{equation*}
\left\{\begin{aligned}
&t_mF(D^2u_m)+(1-t_m)\Delta u_m=0&& ~~\mbox{in}~~B_1;\\
&u_m=g&& ~~\mbox{on}~~\partial B_1.
\end{aligned}\right.
\end{equation*}
Then we have
\begin{equation*}
\left\{\begin{aligned}
&\bar tF(D^2\bar{u})+(1-\bar t)\Delta \bar u=0&& ~~\mbox{in}~~B_1;\\
&\bar u=g&& ~~\mbox{on}~~\partial B_1.
\end{aligned}\right.
\end{equation*}
Hence, $\bar{t}\in A$.

From above arguments, we have shown that $A$ is both open and closed. Therefore, $1\in A$. That is, there exists a solution of \cref{e4.0}.~\qed

\section{Proof of Theorem \ref{th1.1}}\label{S4}
In this section, we prove the interior $C^{2,\alpha}$ regularity for viscosity solutions of fully nonlinear elliptic equations. The proof depends on the following existence of solutions.
\begin{lemma}\label{th4.1}
Let $F$ be uniformly elliptic with $F(0)=0$. Suppose that $g\in C(\partial B_1)$. Then there exists a unique solution $u\in C^{2,\alpha}_{loc}(B_1)\cap C(\bar{B}_1)$ of
\begin{equation*}
\left\{\begin{aligned}
&F(D^2u)=0&& ~~\mbox{in}~~B_1\subset \mathbb{R}^2;\\
&u=g&& ~~\mbox{on}~~\partial B_1.
\end{aligned}\right.
\end{equation*}
Moreover, we have
\begin{equation*}
\|u\|_{C^{2,\alpha}(\bar{B}_{r})}\leq C\|g\|_{L^{\infty}(\partial B_1)},~\forall ~0<r<1,
\end{equation*}
where $C$ depends only on $\lambda,\Lambda$ and $r$.
\end{lemma}
\proof We use \Cref{th3.3} and an approximation argument to prove the lemma. Take two sequences of $F_m\in C^{\infty}$ and $g_m\in C^{\infty}(\partial B_1)$ such that
\begin{itemize}
  \item $F_m(0)=0$, $F_m$ are uniformly elliptic with $\lambda,\Lambda$;
  \item $F_m\to F~~\mbox{ uniformly in any compact set of }~~\mathcal{S}^2$;
  \item $g_m\to g~~\mbox{ uniformly on }~~\partial B_1$;
  \item $g_m$ have the same modulus of continuity on $\partial \Omega$.
\end{itemize}
In fact, above sequences can be constructed by the standard mollification. By \Cref{le4.1}, there exists a sequence of solutions $u_m\in C^{2,\alpha}(\bar{B}_1)$ of
\begin{equation}\label{e4.3}
\left\{\begin{aligned}
&F_m(D^2u_m)=0&& ~~\mbox{in}~~B_1;\\
&u_m=g_m&& ~~\mbox{on}~~\partial B_1.
\end{aligned}\right.
\end{equation}

From \Cref{th3.1} and maximum principle, for any $0<r<1$,
\begin{equation*}
\|u_m\|_{C^{2,\alpha}(\bar{B}_r)}\leq C\|u_m\|_{L^{\infty}(B_1)}\leq C\|g_m\|_{L^{\infty}(\partial B_1)}
\leq C\|g\|_{L^{\infty}(B_1)},
\end{equation*}
where $C$ depends on $r$ and is independent of $m$. Hence, up to a subsequence, there exists $u\in C^{2,\alpha}_{loc}(B_1)$ such that
\begin{equation}\label{e4.4}
u_m\to u\quad\mbox{in}~~C^2(\bar{B}_r),~\forall ~0<r<1.
\end{equation}
In addition, $u_m$ are equicontinuous on $\bar{B}_1$ since $g_m$ have the same modulus of continuity (see \cite[Proposition 4.14]{MR1351007}). Hence,
\begin{equation}\label{e4.5}
u_m\to u\quad\mbox{in}~~C(\bar{B}_1).
\end{equation}

With the aid of \cref{e4.4} and \cref{e4.5}, by taking the limit in \cref{e4.3}, we have
\begin{equation*}
\left\{\begin{aligned}
&F(D^2u)=0&& ~~\mbox{in}~~B_1;\\
&u=g&& ~~\mbox{on}~~\partial B_1.
\end{aligned}\right.
\end{equation*}
That is, $u$ is the desired solution. ~\qed

Now, we give the~\\
\noindent\textbf{Proof \Cref{th1.1}.} Since $u\in C(\bar{B}_1)$, by \Cref{th4.1}, there exists a unique solution $v\in C^{2,\alpha}_{loc}(B_1)\cap C(\bar{B}_1)$ of
\begin{equation*}
\left\{\begin{aligned}
&F(D^2v)=0&& ~~\mbox{in}~~B_1;\\
&v=u&& ~~\mbox{on}~~\partial B_1.
\end{aligned}\right.
\end{equation*}
From the uniqueness of solutions (see \cite[Corollary 5.4]{MR1351007}), $v\equiv u$. That is, $u\in C^{2,\alpha}(B_1)$ and the estimate \cref{e1.1} holds.~\qed

\printbibliography

\end{document}